\documentclass{article}
\usepackage{amscd, amsmath, amssymb}

\title{Birational geometry of quadrics in characteristic 2}
\author{Burt Totaro}
\date{  }

\def\Z{\text{\bf Z}}

\def\C{\text{\bf C}}
\def\P{\text{\bf P}}
\def\F{\text{\bf F}}
\def\arrow{\rightarrow}

\def\qed{\ QED \vspace{.1in}}

\def\deg{\text{deg}}
\def\ql{\text{ql}}

\def\sm{\text{sm}}

\def\dimes{\dim_{\text{es}}}

\def\an{\text{an}}
\begin{document}
\maketitle

\newtheorem{theorem}{Theorem}[section]
\newtheorem{corollary}[theorem]{Corollary}
\newtheorem{lemma}[theorem]{Lemma}
\newtheorem{conjecture}[theorem]{Conjecture}

The theory of quadratic forms can be regarded as studying
an important special case of the general problem
of birational classification of algebraic varieties.
A typical example of a Fano fibration in minimal model theory
is a quadric bundle, which can be viewed without loss
of information as a quadric hypersurface over the function field
of the base variety. This description feeds the problem
into the theory of quadratic forms, which seeks to classify
quadrics over an arbitrary field up to birational
equivalence. One can hope to extend the successes of quadratic form theory
to broader areas of birational geometry.

In the other direction,
birational geometry emphasizes the interest of particular problems
in the theory of quadratic forms. In the study
of Fano fibrations in positive
characteristic, one cannot avoid fibrations whose generic
fiber is nowhere smooth over the function field of the base variety,
such as
``wild conic bundles'' \cite{MS}. That can happen even when
the total space is smooth. In this paper,
we solve two of the main problems of birational geometry
for quadrics which are nowhere smooth over a field.
Such quadrics are called quasilinear; they exist only
in characteristic 2.

First, the ``quadratic Zariski problem'' \cite{Ohm} has a positive answer
for quasilinear quadrics: for any quasilinear quadrics $X$ and $Y$
of the same dimension over a field such that there are rational
maps from $X$ to $Y$ and from $Y$ to $X$, $X$ and $Y$ must be
birational (Theorem \ref{zariski}). Next, a quasilinear quadric
over a field $k$ is ruled, that is, birational to $Z\times \P^1$
for some variety $Z$ over $k$, if and only if its first
Witt index is greater than 1 (Theorem \ref{ruled}). At the same time,
we find that if a quasilinear quadric is ruled, then it
is ruled over a lower-dimensional quadric $Z$. All these
statements are conjectured for arbitrary quadrics, but they are
wide open for other classes of quadrics, for example in
characteristic not 2.

The proofs begin by extending Karpenko and Merkurjev's
theorem on the essential dimension of quadrics to characteristic 2
\cite{KM}. The extension has already been made for smooth quadrics
in the book by Elman, Karpenko, and Merkurjev
\cite{EKM}. Here we cover all quadrics in characteristic 2,
smooth or not. As in other advances on quadratic forms
over the past decade, the proofs use the Chow groups
of algebraic cycles on quadrics and products of quadrics.

Karpenko and Merkurjev's theorem is a major achievement
of the theory of quadratic forms. It includes
Karpenko's theorem that an anisotropic quadric
with first Witt index equal to 1 is not ruled
\cite[Theorem 6.4]{Karpenko},
which we also extend here to characteristic 2.
This is a strong nontriviality statement about the birational
geometry of quadrics.

As an application, we extend Koll\'ar's results on the
birational geometry of conics (including non-smooth conics
in characteristic 2) \cite[section 4]{Kollarcover} to quadrics of any
dimension (Corollary \ref{uniruling}). Finally, besides
proving the conjectured characterization of ruledness
(Conjecture \ref{quadruled}) for all quasilinear quadratic forms,
we check it for non-quasilinear quadratic
forms in characteristic 2 of dimension at most 6, using the
results of Hoffmann and Laghribi (section \ref{questions}).
It is known for quadratic forms of dimension at most 9 in characteristic
not 2 \cite{Totaroauto}.

The inspiration for this work was Hoffmann and Laghribi's
series of papers showing that
the classical theory of quadratic forms (usually restricted
to characteristic not 2, or at most to nonsingular quadratic forms
in characteristic 2)
admits a rich generalization to singular quadratic forms
in characteristic 2,
including the extreme case
of quasilinear forms \cite{Hoffmanndiag,HLTrans,
HL}.

Thanks to Detlev Hoffmann, Bruno Kahn, Nikita Karpenko, and Ahmed Laghribi
for useful conversations.

\section{Notation}

A quadratic form $q$ means a homogeneous polynomial of degree 2 on
a finite-dimensional vector space over a field.
We call $q$ {\it nonsingular }if the associated bilinear
form $b(x,y)=q(x+y)-q(x)-q(y)$ is nondegenerate.
A quadratic form $q$ over a field $k$ 
of characteristic 2 is called {\it quasilinear }(or
totally singular) if the associated bilinear form is identically zero,
or equivalently if $q$ has the form $a_1x_1^2+\cdots+a_nx_n^2$
in coordinates. Equivalently, a nonzero quadratic form $q$ is quasilinear
if and only if the associated projective quadric $Q=\{q=0\}$ becomes
isomorphic to the double hyperplane $x^2=0$ over the algebraic
closure of $k$.
For any quadratic form $q$ on a vector space $V$, the {\it radical }of $q$
(called the quasilinear part $\ql(q)$ in characteristic 2) is the linear
subspace of the elements of $V$ orthogonal to all of $V$. The non-smooth
locus of the associated projective quadric $Q$ over $k$
is the subquadric associated
to the radical of $q$. Two quadratic forms over $k$
are called {\it similar }if one is isomorphic
to a scalar multiple of the other.
A {\it subform }of a quadratic form means any linear subspace with
the restricted quadratic form. This differs from
Hoffmann-Laghribi's terminology, where a form $\alpha$ is called
a subform of $\beta$ if it is an orthogonal summand of
$\beta$ \cite[Definition 2.8]{HL}. 

A quadratic form $q$ over a field is called {\it isotropic }if there
is a nonzero vector $x$ such that $q(x)=0$. The function field
$k(q)$ means the field $k(Q)$ of rational functions on the associated
projective
quadric $Q$. The {\it total index }$i_t(q)$
is defined to be the maximum dimension of an isotropic
subspace (a linear subspace on which $q$ is identically zero).
(For nonsingular forms,
the total index
coincides with the ``Witt index'', the maximum number of orthogonal
copies of the hyperbolic plane that a form contains
\cite[Definition 2.4]{HL}.)
We define the {\it first Witt index }$i_1(q)$
of an anisotropic quadratic form $q$ of dimension at least 2
over a field to be the total index
of $q$ over $k(q)$. We sometimes write $i_1(Q)$ for $i_1(q)$.
We always have $i_1(q)\geq 1$. Section \ref{questions}
helps to show the meaning of the first Witt index for forms
of low dimension in characteristic 2.

Springer's theorem says that an anisotropic quadratic form
over a field $k$ remains anisotropic over any field extension
of odd degree. The usual proof 
\cite[Theorem VII.2.7]{Lam} works
in any characteristic.

We  make some elementary observations which may be helpful. 
An anisotropic quadratic form is automatically nonsingular
in characteristic not 2, but not necessarily in characteristic 2.
For one thing, a nonsingular quadratic form in characteristic 2
has even dimension (because the associated bilinear form is
a nondegenerate alternating form), and so a codimension-1
subform of a nonsingular anisotropic form must have a 1-dimensional
quasilinear part. Moreover, there are plenty of anisotropic forms
with larger quasilinear part over suitable fields: for example,
for any field $k_0$ of characteristic 2, the quasilinear form
$\langle t_1,\ldots,t_n\rangle=t_1x_1^2+\cdots+t_nx_n^2$
is anisotropic over the function field $k_0(t_1,\ldots,t_n)$,
as one easily checks. Over a perfect field of characteristic 2,
an anisotropic quasilinear form has dimension at most 1.
But it would be wrong to restrict oneself to perfect
fields in characteristic 2, since that excludes all function
fields of positive dimension.

\section{Basic properties of quasilinear forms}

In this section we give quick proofs of
some known properties of quasilinear quadratic forms, for later use.
These were proved earlier
by Hoffmann and Laghribi \cite[Proposition 5.3]{Hoffmanndiag},
\cite[Corollaire 3.3]{Laghribi}, \cite[Corollary 2.7(2)]{HL}.

\begin{lemma}
\label{sep}
An anisotropic quasilinear form over a field $k$
of characteristic 2 remains anisotropic
over any separable extension of $k$.
\end{lemma}

{\bf Proof. }It suffices to show 
that an anisotropic quasilinear quadratic form $q$ over a field
$k$ remains anisotropic over a finite Galois extension $l/k$.
The point is that the set of isotropic vectors for $q$ over $l$
is an $l$-linear subspace which is invariant under the Galois group.
Therefore it comes from a $k$-linear subspace. \qed

\begin{lemma}
\label{quasi}
Let $X$ and $Y$ anisotropic quadrics over a field $k$
of characteristic 2 such that
$Y$ becomes isotropic over the function field $k(X)$. If $Y$ is
quasilinear, then $X$ is quasilinear.
\end{lemma}

{\bf Proof. }By Lemma \ref{sep}, since $Y$ has no $k$-points, it
has no rational points over the separable closure of $k$. If $X$
is not quasilinear, then a dense open subset of $X$ is smooth over $k$,
and so the rational points of $X$ over the separable closure of $k$
are dense in $X$. Therefore there is no rational map from $X$ to $Y$
over $k$; that is, $Y$ is anisotropic over $k(X)$. \qed

\begin{lemma}
\label{hyp}
Let $q$ be a anisotropic quadratic form
over a field $k$ of characteristic 2.
If $q$ is not quasilinear, then $q_{k(q)}$ is the orthogonal sum
of $i_1(q)$ hyperbolic planes with an anisotropic form.
\end{lemma}

{\bf Proof. }By Lemma \ref{quasi}, the quasilinear part
$\ql(q)$ remains anisotropic over $k(q)$. By definition
of the first Witt index $i_1(q)$, the maximum dimension
of an isotropic subspace in $q_{k(q)}$ is $i_1(q)$. Since
$\ql(q_{k(q)})$ is anisotropic,
for each nonzero isotropic vector $x$ in $q_{k(q)}$ there is a vector
in the space which is not orthogonal to $x$. This allows us to split
off hyperbolic planes as we want. \qed

\section{Many quadrics are not ruled}

Karpenko showed that an anisotropic quadric with first Witt index 1
over a field of characteristic not 2 is not ruled
\cite[Theorem 6.4]{Karpenko}.
In this section,
we extend Karpenko's theorem
to characteristic 2. 

For any $n$-dimensional variety $X$ over a field $k$,
we define the group  of mod 2
correspondences from $X$ to itself 
to be the Chow group $CH_n(X\times X)/2$ of $n$-dimensional algebraic
cycles on $X\times X$ modulo rational equivalence. For the basic
properties of Chow groups, we refer to Chapter 1 of Fulton's book
\cite{Fulton}.
Note that correspondences are often defined only for
$X$ smooth and proper over $k$, in which case one can
compose correspondences. Nonetheless,
we can use the group of correspondences even for non-smooth quadrics in
characteristic 2. For any correspondence $\alpha\in CH_n(X\times X)/2$, define
the degrees $\deg_1(\alpha)$ and $\deg_2(\alpha)$ 
in $\Z/2$ to be the pushforwards of $\alpha$
via the two projections to $CH_n(X)/2\cong \Z/2$.

Here is our extension of Karpenko's theorem \cite[Theorem 6.4]{Karpenko}
to any characteristic. The proof follows the original one but is
modified in order to avoid the restriction to smooth quadrics. For quasilinear
quadrics, the original proof does not work at all, but there is a simple
direct argument. The extension to {\it smooth }quadrics in
characteristic 2 was proved recently by Elman, Karpenko, and Merkurjev
\cite[Theorem 72.3]{EKM}. We say that a variety $X$ over a field $k$ is {\it
ruled }if it is birational to $Y\times \P^1$ for some variety
$Y$ over $k$.

\begin{theorem}
\label{karp}
Let $X$ be an anisotropic quadric (not necessarily smooth)
of dimension $n$ over a field $k$,
and suppose that the first Witt index $i_1(X)$ is equal to 1. For
any correspondence $\alpha\in CH_n(X\times X)/2$, we have
$\deg_1(\alpha)=\deg_2(\alpha)$ in $\Z/2$.
\end{theorem}

\begin{corollary}
\label{nonruled}
Let $X$ be an anisotropic quadric (not necessarily smooth) over a field $k$.
Suppose that $X$ 
has first Witt index equal to 1. 
Then $X$ is not ruled over $k$.
\end{corollary}

{\bf Proof }(Corollary \ref{nonruled}).
Every rational map
from $X$ to itself determines a correspondence by taking the closure
of the graph. Thus Theorem \ref{karp} implies that every rational map
from $X$ to itself has odd degree. If
$X$ were birational to a product $Y\times \P^1$ over $k$, then projecting
to $Y$ and then including $Y\cong Y\times \{ a \}$
for some point $a\in \P^1(k)$ into $Y\times \P^1$ would give
a rational endomorphism of $X$ with degree 0. (As written, this proof
assumes that $k$ is infinite, which is reasonable since the situation
is trivial for finite fields: every anisotropic quadric
over a finite field
is smooth of dimension at most 0.) \qed

{\bf Proof }(Theorem \ref{karp}).
We first prove this for $k$ of characteristic 2 and $q$ quasilinear.
This case differs from all others in that $X$ is nowhere smooth over $k$
and $X$ does not become rational over any extension field. Nonetheless,
we can argue as follows. Since $q$ has first Witt index equal to 1,
the quadric $X$ over $k(X)$ becomes the projective cone over an
anisotropic quadric $X'$. For any variety $S$ of the same dimension $n$
as $X$ with a morphism $S\arrow X$ of odd degree,
$X'$ remains anisotropic over $k(S)$ by Springer's theorem.
Therefore $X$ has only one rational point over
$k(S)$, the base of the cone. Equivalently, there is only one rational
map from $S$ to $X$ over $k$. Thus, if $S$ is a prime correspondence
(a subvariety of $X\times X$ of dimension $n$) with
$\deg_1(S)=1\pmod{2}$, then the second projection $S\arrow X$ must be
the same map as the first projection. In particular, $\deg_2(S)=1\pmod{2}$
and the theorem is proved (since the same proof shows that $\deg_2(S)$
odd implies $\deg_1(S)$ odd). Clearly it suffices to consider
prime correspondences, as we have done.

We now consider the remaining case, where $k$ can have any characteristic
but $q$ is not quasilinear.
For any correspondence $\alpha\in CH_n(X\times X)/2$,
we will show that 
$\deg_1(\alpha)=\deg_2(\alpha)$ in $\Z/2$. 

The assumption that $i_1(q)=1$ means that over the function field
$k(q)$, the form $q_{k(q)}$ has total index 1. Since $q$ is not quasilinear,
the obvious ``generic'' rational point of $X$ over $k(q)$ is not orthogonal
to the whole vector space $V$ where $q$ is defined. Therefore
$q_{k(q)}$ can be written as an orthogonal sum
$H\perp q'$, where $H$
is the hyperbolic plane. Since $q_{k(q)}$ has total index 1,
the form $q'$ is anisotropic.
From now on, write $X$ for the quadric associated to $q_{k(q)}$
and $X'$ for the codimension-2 subquadric associated to $q'$.

Let $A=X-X'$. It is immediate that $A$ is smooth over $k(q)$.
The trick is to observe that the open variety $A$ behaves like a compact
variety ``modulo 2'', in the sense that we have a well-defined degree map
$$\deg: CH_0(A)/2\arrow \Z/2.$$
To prove that, we use the standard exact sequence for Chow groups
\cite[Proposition 1.8]{Fulton}, 
$$CH_0X'\arrow CH_0X\arrow CH_0A\arrow 0.$$
So it suffices to show that every zero-cycle on $X'$ has even degree.
That is immediate from Springer's theorem, since $X'$ is an anisotropic
quadric. 
Likewise, let us now show that we get well-defined numbers $\deg_i(\alpha)$
for any correspondence $\alpha\in CH_n(A\times A)/2$ and $i=1$ or 2.
Since
$X$ is rational over $k(q)$, $X'$ remains anisotropic over
the function field of $X$. So every $n$-dimensional cycle in $X\times X'$
has even degree over the first factor $X$, and trivially every
$n$-dimensional cycle in $X'\times X$ has zero degree over the first
factor $X$. So the standard exact sequence for Chow groups shows that
$\deg_1(\alpha)$ is defined in $\Z/2$ for any correspondence $\alpha$
in $CH_n(A\times A)$. The same goes for $\deg_2(\alpha)$.

Let $X^{\sm}$ denote the smooth locus of $X$. 
Define a homomorphism $\delta:CH_n(X\times X)/2\arrow \Z/2$
as the composition
$$CH_n(X\times X)/2 \arrow CH_n(X^{\sm}\times X^{\sm})/2
\arrow CH_0(X^{\sm})/2\arrow \Z/2,$$
where the second map is pullback by the diagonal map
$X^{\sm}\arrow X^{\sm}\times X^{\sm}$,
which is a local complete intersection map
since $X^{\sm}$ is smooth over $k(q)$ \cite{Fulton}, and the third map
is our degree map. Since the degree map is well-defined on the
open subset $A$ of $X^{\sm}$, $\delta(\alpha)$ only depends on the
restriction of $\alpha$ to $CH_n(A\times A)/2$.

In coordinates, the quadric $X$ over $k(q)$ is defined
by an equation $x_1x_2+q'=0$, and $A$ is the open subset of $X$
where $x_1$ and $x_2$ are not both zero. Let $A_1$ and $A_2$ be
the hypersurfaces in $A$ defined by $x_1=0$ and $x_2=0$,
respectively. Clearly $A_1$ and $A_2$ are disjoint, and
$A-A_1$ and $A-A_2$ are isomorphic to open subsets of affine space.
So homotopy invariance \cite[Theorem 3.3(a)]{Fulton}
and the standard exact sequence for Chow groups 
imply that
$$CH_n((A-A_1)\times A)=CH_0A$$
and
$$CH_n(A\times (A-A_2))=CH_0A.$$
Since the quadric $X$ has a rational point $p=[1,0,\ldots,0]$ in the smooth
locus, we have $CH_0(X)/2\cong \Z/2
\cdot [p]$
and hence, using our degree map, $CH_0(A)/2\cong \Z/2\cdot [p]$.
 We can then use
the Mayer-Vietoris sequence for the open covering of $A\times A-A_1\times A_2$
by $(A-A_1)\times A$ and $A\times (A-A_2)$ to compute that
$$CH_n(A\times A - A_1\times A_2)/2 \cong \Z/2\oplus \Z/2,$$
where the isomorphism  is given by mapping $\alpha$ to
$(\deg_1(\alpha),\deg_2(\alpha))$.

We now complete the proof. Let $\alpha$ be any correspondence from
$X$ to itself defined over $k$. Then $\delta(\alpha)=0$ in $\Z/2$,
since this is defined by restricting $\alpha$ to a zero-cycle
on the diagonal $X^{\sm}\subset X^{\sm}\times X^{\sm}$
and
every zero-cycle on $X^{\sm}$ over $k$ has even degree (as $X$ is anisotropic
over $k$). Then extend scalars to $k(q)$, where we have
a rational point $p\in A\subset X$, and
consider the restriction of $\alpha$ to $CH_n(A\times A)/2$.
Since $\deg_1(A\times p)=1$, the cycle
$$\widetilde{\alpha}:=\alpha-\deg_1(\alpha)(A\times p)
-\deg_2(\alpha)(p\times A)$$
maps to zero under the above homomorphism to $\Z/2\oplus \Z/2$. Therefore,
$\widetilde{\alpha}$ is the image of some element of $CH_n(A_1\times A_2)$.
Since $A_1$ and $A_2$ are disjoint, this element restricts to zero on the
diagonal of $A$, and so we have $\delta(\widetilde{\alpha})=0$ in $\Z/2$.
On the other hand, it is clear that
$$\delta(\widetilde{\alpha})\equiv 
\delta(\alpha)-\deg_1(\alpha)-\deg_2(\alpha)\pmod{2},$$
and so $\deg_1(\alpha)\equiv \deg_2(\alpha)\pmod{2}$. 
\qed

Let us record a corollary on the birational geometry
of quadrics in any characteristic, slightly stronger
than non-ruledness. The statement here extends Koll\'ar's results on
conics (which specifically included non-smooth conics in characteristic 2)
\cite[section 4]{Kollarcover}.
We say that a {\it uniruling }of an $n$-dimensional variety $X$
over a field $k$
is a dominant rational map $Y\times \P^1\dashrightarrow X$ for some
$(n-1)$-dimensional variety $Y$ over $k$.

\begin{corollary}
\label{uniruling}
An anisotropic quadric $X$ over any field $k$ with first Witt index
equal to 1 has no odd degree uniruling. If in addition $X$ is quasilinear,
then it has no separable uniruling.
By contrast, every non-quasilinear quadric (in any characteristic)
has a degree 2 separable uniruling.
\end{corollary}

{\bf Proof. }
Let $X$ be an anisotropic quasilinear quadric with $i_1(X)=1$
over a field $k$ of characteristic 2. Suppose
we have a separable uniruling $Y\times \P^1\dashrightarrow X$. We know that
$X_{k(X)}$ is the projective cone over an anisotropic quadric $X'$.
The function field $k(Y\times \P^1)$ is a separable extension of $k(X)$,
and so the quasilinear quadric $X'$ remains anisotropic over $k(Y\times
\P^1)$ by Lemma \ref{sep}.
So $X_{k(Y\times \P^1)}$ is the projective
cone over an anisotropic quadric and hence has exactly one rational point.
That is, there is only one rational map from $Y\times \P^1$ to $X$
over $k$, which must be the given dominant map. That is a contradiction, since
we could also project $Y\times \P^1$ to $Y\times \{a\}$ for some
$a\in\P^1(k)$ and then compose with our map to $X$. Thus $X$ is not
separably uniruled. (In particular, $X$ has no odd degree uniruling.)

Next, let $X$ be an anisotropic quadric over any field $k$
such that $i_1(X)=1$ and $X$ is not
quasilinear. Let $n$ be the dimension of $X$ over $k$.
If there is an odd degree separable
uniruling $f:Y\times \P^1\dashrightarrow X$, let $p$ be a general point
in $\P^1(k)$ and let $\pi$ be the projection map from $Y\times \P^1$
to $Y\times \{ p \}$. Consider the rational map
from $Y\times \P^1$ to $X\times X$ given by $(f,f\pi)$ and take the closure
of the image. This gives an $n$-dimensional cycle
 on $X\times X$ with odd degree over the
first factor and zero degree over the second factor, contradicting
Theorem \ref{karp}.
So $X$ has no odd degree uniruling.

Finally,
let $X\subset \P^{n+1}$ be a non-quasilinear quadric over a field $k$ of
any characteristic. We repeat the argument from Koll\'ar
\cite[Remark 4.1.7]{Kollarcover}.
A dense open subset of $X$ is smooth over $k$. A general line
in $\P^{n+1}$ intersects $X$ transversely, and so the intersection points
are in $k$ or in a degree 2 separable extension $E$ of $k$. Thus $X$
is rational over $E$. A fortiori, $X$ has a degree 2 separable uniruling
over $k$.
\qed

\section{The Chow group of zero-cycles on a quadric}

In this section we extend two steps in Karpenko-Merkurjev's
proof to arbitrary quadrics in characteristic 2. First, we show that
the Chow group of zero-cycles
injects into the integers by the degree map.
In characteristic not 2, this is due
to Swan \cite{Swan} and Karpenko \cite[Prop. 2.6]{Karpenkozero}.
We will extend Kahn's proof \cite[Th\'eor\`eme 3.2]{Kahn}, 
a simple variant
of the proof of Springer's theorem, to characteristic 2.
Later in the section, we determine the first Witt index
of ``generic'' subforms of any quadratic form.

\begin{lemma}
\label{chow}
Let $Q$ be any anisotropic quadric over a field $k$,
or any isotropic quadric of dimension greater than 0. Then the degree map gives
an injection from the Chow group of zero-cycles on $Q$ to the integers.
\end{lemma}

The conclusion clearly fails for a smooth isotropic quadric of dimension 0
(which consists of two points).

{\bf Proof. }
This is known in characteristic not 2, and so we assume
that $k$ has characteristic 2 (the proof in characteristic not 2
is in fact one step shorter).
First suppose that $Q\subset \P^{n+1}$
has a $k$-rational point $p$ and $Q$ has dimension greater than 0.
If $Q$ is a union of two hyperplanes over $k$,
then the conclusion is clear (since the two hyperplanes intersect, by
our dimension assumption). So assume that $Q$ is not 
a union of two hyperplanes. Then
$Q$ does not contain all tangent lines to $Q$ at $p$ (if $Q$ is non-smooth
at $p$, we just mean by this that $Q$ does not contain all lines through $p$).
Therefore, by choosing a general
tangent line to $Q$ at $p$, we see that a line $L$ in $\P^{n+1}$ maps
by the pullback map $CH_1\P^{n+1}\arrow CH_0Q$ 
\cite[Proposition 2.6(a)]{Fulton} to $2p$. For any other
$k$-rational point $q$ in $Q$, the line through $p$ and $q$ is either
contained in $Q$, in which case $q$ is rationally
equivalent to $p$, or it meets $Q$ only in $p$ and $q$, in which
case again $q=L|_Q-p=p$ in $CH_0Q$. Any closed point of $Q$ can be
viewed as the pushforward of an $F$-point $q$ of $Q$ for some finite
extension $F/k$ by the proper map $\pi:Q_F\arrow Q$. By what we have
said, $q$ is rationally equivalent to an integer $m$ times 
the point $p$ on $Q_F$,
and so its pushforward is rationally equivalent to $[F:k]m$ times the
point $p$ on $Q$.

If $Q$ has no $k$-rational point, we will show that that
the pullback map
$CH_1\P^{n+1}=\Z\cdot L\arrow CH_0Q$ is surjective. Any closed point
on $Q$ is the pushforward of an $F$-point $q$ of $Q$ for some
finite extension $F/k$ by the proper map $\pi:Q_F\arrow Q$.
We will show that $\pi_*(q)$ is rationally
equivalent to a multiple of $L|_C$ by induction on $[F:k]$.

If $F$ is not generated by one element over $k$,
then $F$ is not a separable extension of $k$, and so we can find
an intermediate field $k\subset E\subset F$ with $[E:F]=2$. Then
the pushforward of $q$ from $Q_F$ to $Q_E$ is an effective zero-cycle
of degree 2, which spans a line in $\P^{n+1}$ over $E$. If this line
is not contained in $Q_E$, then this zero-cycle of degree 2 on $Q_E$ is
the intersection of a line with $Q_E$, and so its pushforward
to $Q$ over $k$ is $[E:k]L|_C$, as we want. On the other hand,
if this line is contained in $Q_E$, then our zero-cycle of degree 2 on 
$Q_E$ is rationally equivalent to 2 times an $E$-point, and so we are
done by our induction on $[F:k]$.

If the field $F$ is generated by one element
over $k$, say
$F=k(\alpha)$, then $q$ lies on a rational curve $C$
(the image of a map $\P^1\arrow \P^{n+1}$ over $k$)
of degree $d$ which is less than $[F:k]$. This curve is found by writing
the coordinates in $F$ as polynomials in $\alpha$. This curve is
not contained in $Q$; otherwise $Q$ would have a $k$-point. So
the pullback of $C$ to $Q$ is on the one hand rationally equivalent
to $dL|_C$ and on the other is equal to $\pi_*(q)+z$, where $z$
is an effective zero-cycle of degree $2d-[F:k]<[F:k]$. Thus,
by our induction on $[F:k]$, $z$ is rationally equivalent to a multiple
of $L|_C$, and so the same goes for $\pi_*(q)$.
\qed

We now extend another step in Karpenko-Merkurjev's
argument \cite[Proposition 1.2]{KM}
to arbitrary quadrics in characteristic 2: the determination
of the first Witt index of ``generic'' subforms of a given
form. These subforms are defined over a purely transcendental
extension of the ground field.

\begin{lemma}
\label{generic}
Let $q$ be an anisotropic quadratic form over a field $k$,
and let $j$ be an integer with $0\leq j\leq \dim q-2$. Then
there is a purely transcendental extension $F/k$ and
a codimension-$j$ subform $s$ of $q_F$ such that
$$i_1(s)=\begin{cases} i_1(q)-j & \text{if }i_1(q)>j; \\
    1 & \text{if }i_1(q)\leq j.
\end{cases}$$
\end{lemma}

{\bf Proof. }
We assume that $k$ has characteristic 2, the result being
known for characteristic not 2.
It suffices to prove the lemma for $j=1$.
Let $F$ be the function field of the projective space of
hyperplanes in the vector space of $q$; clearly $F$ is a purely
transcendental extension of $k$. We have a canonical
codimension-1 subform $s$ of $q_F$ (the ``generic'' subform).
We want to compute the total index (maximal dimension of an
isotropic subspace) of $s$ over the function field $E/F$
of the corresponding projective quadric $S$. Clearly $E$
is the function field over $k$ of the variety of pairs $(p,H)$ where
$p\in Q$ and $H$ is a hyperplane in $\P^{n+1}$ that contains $p$,
where $Q\subset \P^{n+1}$
is the projective quadric associated to $q$. So (by choosing
$p$ first and then $H$, so to speak) we can also view
$E$ as the function field over $k(q)$ of the variety of hyperplanes
in $\P^{n+1}$ that contain the canonical (``generic'') point $p_0$ in $Q(k(q))$.
In particular, $E$ is a purely transcendental extension of $k(q)$.
From this point of view, it becomes easy to compute the total
index of $s$ over the field $E$. Clearly the quadric $S_E$
is a codimension-1 subquadric of $Q_E$ that contains
the rational point $p_0$, and so $s_E$ has total index at least 1.

The form $q_{k(q)}$ has total index equal to $i_1(q)$,
by definition. Since $E$ is a purely transcendental extension
of $k(q)$, $q_E$ also has total index equal to $i_1(q)$. If
$i_1(q)$ is equal to 1, then $s_E$ has total index at most 1
just because it is a subform of $q_E$. Hence $s_E$
has total index equal to 1,
and we are done.

It remains to consider the case where $r:=i_1(q)$ is at least 2.
Clearly the total index
of $s_E$ is at least $r-1$ (as $s_E$ is a codimension-1
subform of $q_E$), and we want to show that equality holds; here
$s_E$ is the ``generic'' codimension-1 subform of $q_{k(q)}$
that contains the isotropic vector $p_0$. It suffices
to find one codimension-1 subform of $q$ over $k(q)$, or 
a purely transcendental extension thereof, which contains
$p_0$ and has total index $r-1$ (rather than $r$). That in turn
follows if we can find one codimension-$r$ subform of $q$ over $k(q)$,
or a purely transcendental extension thereof, which is anisotropic.
This is easy for $q$ quasilinear;
in that case, the maximal isotropic subspace of $q_{k(q)}$ (of dimension
$r$) is unique,
and we can take the desired codimension-$r$ subform to be any complementary
linear subspace over $k(q)$ to the maximal isotropic subspace.

So we can assume $q$ is not quasilinear. Then $q_{k(q)}$ is the
orthogonal sum of $r$ copies of the hyperbolic plane together
with an anisotropic form $\varphi$ by Lemma \ref{hyp}.
Therefore, over the purely transcendental extension $k(q)(t_1,\ldots,t_r)$
of $k(q)$,
$q$ contains a codimension-$r$ subform which is the orthogonal
sum of $\varphi$ and the quasilinear form $\langle t_1,\ldots t_r\rangle$.
One can check by hand that this
subform is anisotropic. That is what we needed. 
\qed

\section{Essential dimension of quadrics}

We now have all the tools needed to extend Karpenko and Merkurjev's
theorem on the essential dimension
to arbitrary quadrics in characteristic 2. This has been
done by Elman, Karpenko, and Merkurjev for {\it smooth }quadrics
in characteristic 2 \cite[Theorem 73.1]{EKM}.

For any anisotropic projective quadric $X$ over a field $k$,
define the {\it essential dimension }$\dimes(X)$
to be $\dim(X)-(i_1(X)-1)$.

\begin{theorem}
\label{ess}
Let $X$ be an anisotropic projective quadric (not necessarily
smooth) over any field $k$,
and let $Y$ be a proper variety over $k$ with all closed points
of even degree. Suppose that $Y$ has a closed point of odd degree
over $k(X)$. Then

(1) $\dimes(X)\leq \dim(Y)$;

(2) if, moreover, $\dimes(X)=\dim(Y)$, then $X$ is isotropic
over $k(Y)$.
\end{theorem}

{\bf Proof. }
For $k$ of characteristic not 2, this is Theorem 3.1 in Karpenko-Merkurjev
\cite{KM}. In characteristic 2, the same proof works,
using the following new tools: 
the calculation of the Chow group of zero-cycles on a quadric
(Lemma \ref{chow}), the calculation of the first Witt index
of a generic subquadric (Lemma \ref{generic}, which extends
Proposition 1.2 in \cite{KM}), and Theorem \ref{karp} on
correspondences of a quadric with first Witt index 1, which
extends Karpenko's theorem (Theorem 3.2 in \cite{KM}).

In more detail, here is the proof.
The assumption means that
there is a correspondence $\alpha:X\rightsquigarrow Y$ of odd degree
over $X$. We can assume that $\alpha$ is a prime correspondence, that is,
an irreducible subvariety $Z$ of $X\times Y$ of dimension $\dim(X)$.
By Springer's theorem, statement (2) follows if we can find
a correspondence $Y\rightsquigarrow X$ of odd degree over $Y$. Also,
it will suffice to show that $X$ is isotropic over $k(Y)$ after
a purely transcendental extension of the field $k$.

Assume first that $i_1(X)=1$, so that $\dimes(X)=\dim(X)$. We prove
(1) and (2) simultaneously by induction on $n=\dim(X)+\dim(Y)$.

For $n=0$, $X$ and $Y$ are varieties of dimension zero over $k$,
with $X$ of degree 2 and $Y$ of even degree over $k$. We are given
another 0-dimensional variety $Z$ of odd degree over $X$ with a map
to $Y$ over $k$. Since $Z$ has degree 2 times an odd number over $k$,
it has odd degree over $Y$. Thus $Z$
is a correspondence from $Y$ to $X$ of odd degree over $Y$,
and the theorem is proved.

So let $n>0$ and assume the theorem for smaller values of $n$.
We first prove (2). Thus, we assume that $\dim(Y)=\dim(X)>0$. It suffices
to show that $\alpha$ has odd degree over $Y$. Suppose that $\alpha$
has even degree over $Y$. By adding a suitable multiple of $x\times Y$
to $\alpha$, where $x$ is a closed point of degree 2 on $X$, we can assume
that $\alpha$ has degree zero over $Y$. Since $CH_0(X_{k(Y)})$ injects
into the integers by the degree map (Lemma \ref{chow}), there is a dense open
subset $U\subset Y$ such that the restriction of $\alpha$ to
$X\times U$ is rationally equivalent to zero. By the standard exact
sequence of Chow groups, $\alpha$ is rationally equivalent to
a correspondence $\alpha':X\rightsquigarrow Y'$ for some lower-dimensional
closed subset $Y'$ of $Y$, still of odd degree over $X$. This
contradicts statement (1) (for a lower value of $n=\dim(X)+\dim(Y)$).

We now prove (1). Thus, assume that $\dim(Y)<\dim(X)$. We will
obtain a contradiction. Replacing $Y$ by the image of the
given prime correspondence $\alpha$ in $Y$,
we can assume that $Z$ maps onto $Y$.

By Lemma \ref{generic}, after replacing $k$ by a purely transcendental
extension, $X$ has a subquadric $X'$ of the same dimension as $Y$
such that $i_1(X')=1$. The hypotheses on $X$ and $Y$ still hold
over the new field $k$. We can pull back $\alpha$ to $X'\times Y$
since $X'$ is a local complete intersection subscheme of $X$
\cite[section 6.2]{Fulton},
yielding a correspondence $X'\rightsquigarrow Y$ of odd degree over $X'$.
By induction, property (2) holds for $X'$ and $Y$; that is, there is
a rational map from $Y$ to $X'$ over $k$, and hence a rational
map $\beta:Y\dashrightarrow X$. We can compose the dominant map $Z\arrow Y$
with the rational map $\beta$ to get a rational map $Z\dashrightarrow X$.
Thus $Z$ gives a correspondence from $X$ to itself, of odd degree over
the first factor and of zero degree over the second factor (because it
maps into the proper subquadric $X'$ of $X$).
But this contradicts Theorem \ref{karp}.

The induction is complete: Theorem \ref{ess} is proved for $i_1(X)=1$.
Now let $i_1(X)$ be arbitrary. After replacing $k$ by a purely
transcendental extension, which does not change the properties
assumed of $X$ and $Y$, $X$ has a subquadric $X'$ of dimension
equal to $\dimes(X)$ such that $i_1(X')=1$, by Lemma
\ref{generic}. The correspondence $\alpha:X\rightsquigarrow Y$ pulls back
to a correspondence from $X'$ to $Y$, still of odd degree over $X'$.
By the first part of the proof, $\dimes(X)=\dim(X')\leq \dim(Y)$. 
If $\dim(X')=\dim(Y)$, then the first part of the proof shows that
$X'$ and hence $X$ have a rational point over $k(Y)$. 
\qed

In particular, an anisotropic quadric $X$ cannot be ``compressed''
below $\dimes(X)$, in the sense that there is no rational map from
$X$ to any variety of dimension less than $\dimes(X)$ with
all closed points of even degree.

We now consider the most natural situation, where the variety $Y$
is also a quadric.

\begin{theorem}
\label{essquad}
Let $X$ and $Y$ be anisotropic quadrics (not necessarily smooth)
over any field $k$. Suppose that $Y$ is isotropic over $k(X)$.
Then

(1) $\dimes(X)\leq \dimes(Y)$;

(2) moreover, the equality $\dimes(X)=\dimes(Y)$ holds if and only
if $X$ is isotropic over $k(Y)$.
\end{theorem}

{\bf Proof. }
For $k$ of characteristic not 2, this is Theorem 4.1 in
Karpenko-Merkurjev \cite{KM}. Their proof
works in characteristic 2, as follows.
Let $Y'$ be a subquadric
of $Y$ with $\dim(Y')=\dimes(Y)$. 
The definition of $i_1(Y)$ implies that $Y'$
becomes isotropic over the function field $k(Y)$.
For any anisotropic quadrics $A,B,C$ over any field $k$
with $A$ isotropic over $k(B)$ and $B$ isotropic over $k(C)$,
$A$ is isotropic over $k(C)$ \cite[Proposition 22.16]{EKM}.
Since $Y$ is isotropic over $k(X)$, $Y'$ is isotropic over $k(X)$.
By Theorem \ref{ess}, $\dimes(X)\leq \dim(Y')$,
which proves statement (1).

If $\dimes(X)=\dimes(Y)$, then Theorem \ref{ess} implies that
$X$ is isotropic over $k(Y')$. Since $Y'$ is isotropic over $k(Y)$,
$X$ is isotropic over $k(Y)$.
Conversely, if $X$ is isotropic
over $k(Y)$, then applying statement (1) with $X$ and $Y$ reversed
shows that $\dimes(Y)\leq \dimes(X)$, and hence the equality holds.
Thus (2) is proved.
\qed

Laghribi pointed out the following application (Corollary \ref{app}), of which
the last part is new for quasilinear forms. It was 
proved by Izhboldin in characteristic
not 2 \cite[Theorem 0.2]{Izhboldinmot}, as part of his construction
of fields with $u$-invariant 9. Later Corollary \ref{app}
in characteristic not 2 became
an easy consequence of the Karpenko-Merkurjev theorem.
It was proved by Hoffmann-Laghribi
for non-quasilinear forms in characteristic 2 \cite[Theorem 1.3]{HL},
but Theorem \ref{essquad} allows an easier approach.

\begin{corollary}
\label{app}
Let $q$ and $r$ be anisotropic quadratic forms over any
field $k$ such that $\dim r>2^n$ for some integer $n$.
If $q$ has dimension at most $2^n$,
then $q$ is not isotropic over $k(r)$.
If $q$ has dimension $2^n+1$,
then $q$ is isotropic over $k(r)$ if and only if $r$
is isotropic over $k(q)$.
\end{corollary}

{\bf Proof. }This is immediate from Theorem \ref{essquad} together with
Hoffmann-Laghribi's theorem that $i_1(r)\leq
\dim(r)-2^n$ \cite[Lemma 4.1]{HL}. \qed

\section{Quasilinear forms and the quadratic Zariski problem}

Karpenko's theorem, extended in this paper to characteristic 2 (Corollary
\ref{nonruled}), asserts that every anisotropic quadratic form
with first Witt index equal to 1 is not ruled (that is,
the corresponding projective quadric is not ruled). We conjecture
the converse, in any characteristic (Conjecture \ref{quadruled}).
In this section, we prove the conjecture for quasilinear forms
of any dimension
(Theorem \ref{ruled}).
At the same time, we solve the quadratic Zariski problem for
quasilinear quadrics \cite{Ohm}. That is, two anisotropic quasilinear
quadrics $X$ and $Y$ of the same dimension
such that $X$ is isotropic over
$k(Y)$ and $Y$ is isotropic over $k(X)$ must be
birational (Theorem \ref{zariski}).
In this respect, birational geometry is understood better
for quasilinear quadrics than for any other class of quadrics,
for example those of characteristic not 2. The section concludes
by describing which quasilinear quadrics are regular schemes.

To see that birational geometry is nontrivial for
quasilinear quadrics, note that two quasilinear quadratic forms
whose associated quadrics are birational need not be similar.
We can imitate Ahmad-Ohm's example in characteristic not 2
\cite[Example 1.5]{AO}. Let $k_0$ be any field of characteristic 2,
let $k=k_0(a,b,c)$, and consider the quasilinear forms
$\langle 1,a,b,ab,c\rangle$ and $\langle 1,a,c,ac,b\rangle$ over $k$.
They are not similar, as one checks by hand (write down the properties
that a similarity factor in $k$ must satisfy, and use that the
quasi-Pfister form $\langle \langle a,b,c\rangle\rangle:=
\langle 1,a,b,c,ab,ac,bc,abc\rangle$
is anisotropic over $k$). But the
3-dimensional
quadrics associated to these two forms
are birational over $k$ by Theorem \ref{zariski}, using that
the two forms are neighbors of the same quasi-Pfister form
$\langle \langle a,b,c\rangle\rangle$. (See section
\ref{questions} for the definitions.)

Here is the proposed converse to Karpenko's theorem.

\begin{conjecture}
\label{quadruled}
Let $q$ be an anisotropic quadratic form over a field $k$.
If the first Witt index of $q$ is greater than 1, then the associated
projective quadric is ruled over $k$.
\end{conjecture}

We begin with a special case of the quadratic Zariski problem (Lemma
\ref{one}).
The following lemmas are not stated in the strongest possible form:
for instance, the assumption that $Y$ is isotropic over $k(X)$ can
be omitted from Lemma \ref{one}, in view of Theorem \ref{essquad}.
The lemmas merely serve the purpose of proving Theorems \ref{ruled}
and \ref{zariski}.

\begin{lemma}
\label{one}
Let $X$ and $Y$ be anisotropic quasilinear quadrics over a field $k$
with $X$ isotropic over $k(Y)$ and $Y$ isotropic over $k(X)$. Suppose that
$X$ and $Y$ have the same dimension and that both have first Witt index
equal to 1. Then $X$ and $Y$ are birational over $k$.
\end{lemma}

{\bf Proof. }Let $q$ be a quadratic form whose associated quadric
is $X$. Since $q$ is quasilinear and $i_1(q)=1$, the set of 
isotropic vectors for $q$ over $k(X)$ is a 1-dimensional vector
space over $k(X)$. Equivalently, there is a unique rational map
from $X$ to itself over $k$, which must be the identity map.
The same goes for $Y$. We are assuming that there are rational
maps $f:X\dashrightarrow Y$ and $g:Y\dashrightarrow X$ over $k$. Both maps
are dominant by Theorem \ref{ess} since $X$ and $Y$ have
first Witt index 1. Therefore we can compose them. We must have
$fg=1_Y$ and $gf=1_X$. \qed

\begin{lemma}
\label{same}
Let $X$ be an anisotropic quadric over a field $k$.
Let $Y$ be a subquadric of codimension $i_1(X)-1$.
Let $A$ be any anisotropic
quadric over $k$. Then $A$ has the same
total index over $k(X)$ as over $k(Y)$.
\end{lemma}

{\bf Proof. }By definition of the first Witt index, $Y$ is isotropic
over $k(X)$. That is, there is a rational map $X\dashrightarrow Y$
over $k$. The map is dominant by Theorem \ref{ess}. Using this map,
we can view $k(X)$ as an extension field of $k(Y)$, and so the inequality
$i_t(A_{k(Y)})\leq i_t(A_{k(X)})$ is clear.

For the opposite inequality,
it suffices to show that $i_t(A_{k(Z)})\geq i_t(A_{k(X)})$ for
every codimension-1 subquadric $Z$ in $X$. Since $X$ is anisotropic, so is $Z$.
Therefore $Z$ is irreducible and reduced, and so the scheme $Z$
is regular at its generic
point. Since $Z$ is a Cartier divisor in $X$, it follows that $X$ is regular
at the generic point of $Z$. In other words, the local ring of $X$
at the generic point of $Z$ is a discrete valuation ring.
By the valuative criterion of properness, any rational map from $X$
to a proper variety $B$ over $k$ is defined at the generic point of $Z$.
Let $B$ be the Grassmannian of isotropic linear subspaces of
dimension $i_t(A_{k(X)})$ in a quadratic form associated to $A$.
Then there is a rational map
from $X$ to $B$ over $k$, and hence a rational map from $Z$ to $B$ over $k$.
That is, $i_t(A_{k(Z)})\geq i_t(A_{k(X)})$. \qed

The following statement proves Conjecture \ref{quadruled} for quasilinear
quadrics.

\begin{theorem}
\label{ruled}
Let $X$ be an anisotropic quasilinear quadric over a field $k$
with first Witt index $r$.
Let $Y$ be a subquadric of codimension $r-1$. Then $Y$
has first Witt index 1, and
$X$ is birational to $Y\times \P^{r-1}$ over $k$.
\end{theorem}

{\bf Proof. }Clearly $X$ is isotropic over $k(Y)$.
By definition of the first Witt index, $Y$ is isotropic
over $k(X)$. So there is a rational map $\pi:X\dashrightarrow Y$
over $k$, which is dominant by Theorem \ref{ess}.
We have $i_1(Y)=1$ by Theorem \ref{essquad}.

Apply Lemma \ref{same} to $A=X$. We find
that $X$ has the same total index $r$ over $k(Y)$ as over $k(X)$.
Let $q$ be a quadratic form whose associated quadric is $X$. Since
$q$ is quasilinear, the set of isotropic vectors for $q$ over $k(Y)$
is a vector space over $k(Y)$, which we know 
has dimension $r$. Let $s_1,\ldots,s_r$
be a basis; we can view $s_1,\ldots,s_r$
 as rational maps from $Y$ to the affine
quadric $\{q=0\}$.
Define a rational map $\varphi$ from $Y\times \P^{r-1}$ to $X$ over $k$
by $(y,[a_1,\ldots,a_r])\mapsto [a_1s_1(y)+\ldots +a_rs_r(y)]$.

Think of $k(Y)$ as a subfield of $k(X)$ via the map $\pi$. 
Since $X$ has the same
total index over $k(Y)$ as over $k(X)$, the ($r$-dimensional)
vector space of isotropic
vectors for the quadratic form $q$ over $k(X)$ is the vector space
of isotropic vectors for $q$ over $k(Y)$ tensored over $k(Y)$ with $k(X)$.
The identity map on $X$ corresponds to some element of the first
vector space. Therefore there are rational functions $f_i$ on $X$
such that $[f_1(s_1\circ \pi)+\cdots + f_r (s_r\circ \pi)]$
is the identity map from $X$ to itself. Define a rational
map $\psi$ from $X$ to $Y\times \P^{r-1}$ by $\psi(x)=
(\pi(x),[f_1(x),\ldots,f_r(x)])$; then $\varphi\psi$ is the identity
on $X$. So $\psi$ has degree 1. That is, $\psi$ is birational. \qed

The following statement solves the quadratic Zariski problem for
quasilinear quadrics.

\begin{theorem}
\label{zariski}
Let $X$ and $Y$ be anisotropic quasilinear quadrics over a field $k$
such that $X$ is isotropic over $k(Y)$ and $Y$ is isotropic over $k(X)$.
Then there is a subquadric $X'$ of $X$ such that $X$ is birational
to $X'\times \P^{i_1(X)-1}$
and $Y$ is birational to $X'\times \P^{i_1(Y)-1}$ over $k$.
\end{theorem}

{\bf Proof. }Let $r=i_1(X)$. Let $X'$ be any subquadric of codimension
$r-1$ in $X$. By Theorem \ref{ruled}, $i_1(X')=1$ and 
$X$ is birational to $X'\times
\P^{r-1}$. Likewise, let $s=i_1(Y)$ and let $Y'$ be any subquadric
of codimension $s-1$ in $Y$; then $i_1(Y')=1$ and
$Y$ is birational to $Y'\times \P^{s-1}$.
By Theorem \ref{essquad}, $X'$ and $Y'$ have the same dimension.
By Lemma \ref{one}, $X'$ and $Y'$ are birational. \qed

The proof of Lemma \ref{same}
suggests the question of which quasilinear quadrics
are regular schemes. Such quadrics arise geometrically as the generic
fiber of a nowhere smooth quadric fibration with smooth total space.
The following lemma gives a complete answer,
essentially due to Hoffmann:
a quasilinear quadric is a regular scheme if and only the coefficients
are generic in a strong sense, or again if and only if the splitting
pattern is the generic one. In particular, if a quasilinear quadric
is a regular scheme, then it
is anisotropic and has first Witt index 1; but the converse is false,
as shown by the form $\langle 1,a,b,ab,c\rangle$ over $k=k_0(a,b,c)$.

We define the splitting pattern of a quasilinear quadratic form $q$ over 
a field $k$ as follows.
Let $k_0=k$,
$q_0=q_{\an}$ (the anisotropic part of $q$), $k_{j+1}=k_j(q_j)$,
and $q_{j+1}=((q_j)_{k_{j+1}})_{\an}$, where we stop with $q_h$ such that
$\dim(q_h)\leq 1$. Then the splitting pattern of $q$ is defined
to be the sequence $(\dim(q_0),\ldots,\dim(q_h))$.
We refer to section \ref{questions} for the definition of quasi-Pfister
forms.

\begin{lemma}
Let $q$ be a nonzero 
quasilinear quadratic form over a field $k$. After scaling $q$,
we can assume that it has the form $q=a_1x_1^2+\cdots+a_{n-1}x_{n-1}^2
+x_n^2$. The following are equivalent.

(1) The projective quadric associated to $q$ is a regular scheme.

(2) $da_1,\ldots,da_{n-1}$ are linearly independent over $k$
in the module of differentials $\Omega^1_k$ ($=\Omega^1_{k/\F_2}$).

(3) $da_1\wedge\cdots\wedge da_{n-1}\neq 0$ in $\Omega^{n-1}_k$.

(4) The quasi-Pfister form $\langle\langle a_1,\ldots,a_{n-1}
\rangle\rangle$ is anisotropic over $k$; that is, $1$, $a_1,\ldots,
a_{n-1}$, $a_1a_2,\ldots,a_{n-2}\, a_{n-1},\ldots,
a_1\cdots a_{n-1}$ in $k$ are linearly  independent
over $k^2$.

(5) The form $q$ is anisotropic and has splitting pattern $(n,n-1,
\ldots,1)$.
\end{lemma}

{\bf Proof. }Properties (2) and (3) are equivalent  by the definition
$\Omega^{n-1}_k=\Lambda^{n-1}\Omega^1_k$ (exterior power over $k$).
Hoffmann showed the equivalence of (3), (4), and (5)
\cite[Lemma 8.1(ii), Theorem 7.25]{Hoffmanndiag}.

It remains to show the equivalence of (1) and (2). Since the
quadric $Q$ is a hypersurface in the regular scheme $\P^{n-1}_k$,
the non-regular locus of $Q$ is defined by the equations $q=0$
and $dq=0$, where $dq$ lies in $\Omega^1_{k[x_1,\ldots,x_n]}$.
Since $q=a_1x_1^2+\ldots, a_{n-1}x_{n-1}^2+x_n^2$, we have
$dq=(da_1)x_1^2+\ldots +(da_{n-1})x_{n-1}^2$, which lies in
the subgroup $\Omega^1_k\otimes_k k[x_1,\ldots,x_n]$. Thus
the subscheme of $\P^{n-1}_k$ defined by $dq=0$ is the intersection of
a family of quasilinear quadrics. We read off that
the subscheme defined by $q=0$ and $dq=0$ is empty
if and only if $da_1,\ldots,da_{n-1}$ are linearly
independent in $\Omega^1_k$. \qed

\section{Ruled quadrics in characteristic 2}
\label{questions}

Conjecture \ref{quadruled} proposes that every anisotropic quadratic form
with first Witt index greater than 1 is ruled (that is,
the corresponding projective quadric is ruled). The conjecture is true
for all quasilinear forms (Theorem \ref{ruled}) and for
all quadratic forms of dimension at most 9
in characteristic not 2 \cite{Totaroauto}. In this section, we check
Conjecture \ref{quadruled} for all quadratic
forms of dimension at most 6 in characteristic 2, using the results
of Hoffmann and Laghribi.

Throughout this section we will assume that the base field
$k$ has characteristic 2.

Some evidence for Conjecture \ref{quadruled} is that
a ``stabilized'' version is true (as in characteristic not 2).
We can restrict here to non-quasilinear forms.
Namely, if $Q$ is a non-quasilinear anisotropic quadric
with $r:=i_1(Q)$ greater than 1, let $Q'$ be any
subquadric of codimension
$r-1$. The definition of $i_1(Q)$ implies that $Q'$ is isotropic
over $k(Q)$; that is, there is a rational map $\pi$ from $Q$ to $Q'$
over $k$. The quadric $Q'$ is not quasilinear, for example by
Lemma \ref{quasi}. The map $\pi$ is dominant
by Theorem \ref{ess}. So $Q'$ has a smooth point over $k(Q)$
and hence is rational over $k(Q)$. Likewise,
there is also a
rational map from $Q'$ to $Q$ (the inclusion) whose image meets
the smooth locus of $Q$, and so $Q$ is rational over $k(Q')$.
Thus the product $Q\times
Q'$ is birational both 
to $Q\times \P^a$ and to $Q'\times \P^{r-1+a}$ for some $a$.
So it seems plausible (by analogy with the quadratic Zariski
problem \cite{Ohm}) that $Q$ should be birational to
$Q'\times \P^{r-1}$ (where $r-1\geq 1$) and hence that $Q$ should be ruled.

In order to prove Conjecture \ref{quadruled} in low dimensions,
we will use Lemma \ref{pfister} below, a sufficient condition
for ruledness. The lemma is a characteristic 2 version of a result
by Ahmad and Ohm \cite[after 1.3]{AO}. In order to formulate the lemma,
we need to define the most important classes of quadratic forms
in characteristic 2.

An $n$-fold
Pfister form over a field $k$ of characteristic 2 is the nonsingular
quadratic form defined as
the tensor product
 $$\langle \langle a_1,\ldots,a_{n-1},c]]=
\langle 1,a_1\rangle_b \otimes\cdots
\otimes \langle 1,a_{n-1}\rangle_b \otimes [1,c]$$
for some
$a_i$ in $k^*$ and $c$ in $k$ \cite[section 2.5]{HL}.
Here $\langle 1,a\rangle_b$ denotes the bilinear form
$x_1y_1+ax_2y_2$ and $[1,c]$ is the quadratic
form $x_1^2+x_1x_2+cx_2^2$. Explicitly, for any quadratic form $q$,
the tensor product $\langle a_1,\ldots,a_n\rangle_b\otimes q$
is the quadratic form
$a_1q\perp\cdots\perp a_nq$. An $n$-fold quasi-Pfister form is the
quasilinear quadratic form defined as the tensor product
$$\langle\langle a_1,\ldots,a_n\rangle\rangle=\langle 1,a_1\rangle_b\otimes
\cdots\otimes \langle 1,a_n\rangle_b\otimes \langle 1\rangle$$
for some $a_i$ in $k$ \cite[Definition 4.4]{Hoffmanndiag}.
A {\it Pfister neighbor }means
a quadratic form $q$ similar to a subform of a Pfister form $\varphi$
with $\dim(q)>\dim(\varphi)/2$. {\it Quasi-Pfister neighbors }are defined
analogously.

\begin{lemma}
\label{pfister}
Let $P$ be a Pfister form or
quasi-Pfister form over a field $k$, $P_1$ a nonzero subform of $P$,
and $b_1,\ldots, b_s$ any elements of $k^*$. Then the projective quadric
associated to $b_1P\perp \cdots\perp b_{s-1}P\perp b_sP_1$ is birational
to the product of a projective space with the quadric
$b_1P\perp \cdots\perp b_{s-1}P\perp \langle b_s \rangle$;
in particular, it is ruled if $P_1$ has dimension at least 2.
\end{lemma}

{\bf Proof. }
A Pfister form or quasi-Pfister form $P$ is strongly multiplicative: there is
a ``multiplication''
$xy$ on the vector space of $P$ which is a rational function of $x$ and
linear in $y$ such that $P(xy)=P(x)P(y)$ \cite[Theorem X.2.11]{Lam}.
(For quasi-Pfister forms, the ``multiplication'' can even be taken
to be bilinear, as Albert observed long ago \cite[section 8]{Albert}.)
Given a general point
on the first quadric $b_1P(x_1)+\cdots+b_{s-1}P(x_{s-1})+b_sP(x_s)=0$
(where $x_s\in P_1$), we map it to a point on the second quadric
by noting that $b_1P(x_sx_1)+\cdots +b_{s-1}P(x_sx_{s-1})+b_sP(x_s)^2=0$.
The general fibers of this rational map are linear spaces. \qed

In particular,
any quadratic form divisible by a binary quadratic form is ruled,
and any {\it special Pfister neighbor }or
{\it special quasi-Pfister neighbor }(that is,
a form $b_1P\perp b_2P_1$ where $P$ is a Pfister form
or quasi-Pfister form
and $P_1$ is a nonzero subform of $P$)
is ruled if its dimension is not of the form
$2^a+1$. More generally, Theorem \ref{ruled} implies that every
quasi-Pfister neighbor is ruled if its dimension is not of the form
$2^a+1$, but Lemma \ref{pfister} could be considered more explicit
when it applies.

We now show that an anisotropic form $q$ in characteristic 2
with $i_1(q)>1$ must
be ruled when $q$ has dimension at most 6.
Hoffmann-Laghribi showed that an anisotropic quadratic form $q$
of dimension $2^a+m$ over any field, where $1\leq m\leq 2^a$,
has $i_1(q)\leq m$ \cite[Lemma 4.1]{HL}.
In particular, a 3-dimensional anisotropic form
$q$ has $i_1(q)=1$ and so there is nothing to prove.

So let $q$ be a 4-dimensional anisotropic form with $i_1(q)>1$, hence
$i_1(q)=2$ by Hoffmann-Laghribi. We can apply the following result
of Hoffmann and Laghribi
\cite[Theorem 1.2(1)]{HL}, \cite[Corollary 7.29, 
Remark 7.30]{Hoffmanndiag}. 

\begin{theorem}
\label{neighbor}
Let $q$ be an anisotropic form over a field $k$ of characteristic 2
such that $2^{n+1}-2\leq \dim(q) \leq 2^{n+1}$. If $i_1(q)=\dim(q)-2^n$
(the largest possible), then $q$ is either a Pfister neighbor
or a quasi-Pfister neighbor.
\end{theorem}

In particular, a 4-dimensional anisotropic form $q$ with $i_1(q)>1$
is similar to a Pfister or quasi-Pfister form. In both cases, $q$
is ruled by Lemma \ref{pfister}.

For any 5-dimensional anisotropic form $q$, we have $i_1(q)=1$ and so
there is nothing to prove. So let $q$ be a 6-dimensional form
with $i_1(q)$ greater than 1. Then $i_1(q)=2$, and $q$ is a Pfister
neighbor or quasi-Pfister neighbor by Theorem \ref{neighbor}. If $q$ is
quasilinear, then it is ruled by Theorem \ref{ruled}, and so we
can assume that $q$ is a Pfister neighbor of dimension 6.
The radical of $q$
has dimension 0 or 2. Laghribi described the forms of both types
\cite[Proposition 3.2 and its proof]{Laghribi}. Namely, if $q$ is nonsingular,
then $q$ can be written as
$\langle u,v,w\rangle_b \otimes [1,\Delta]$ for some elements
$u,v,w$ in $k^*$ and $\Delta $ in $k$, while if
the radical of $q$ has dimension 2, then $q$ can be written as
$(r\langle u,v\rangle_b\otimes [1,s])\perp \langle u,v\rangle$.
In both cases, we read off that $q$ is a special Pfister neighbor,
and so $q$ is ruled.

For $q$ of dimension 7, Conjecture \ref{quadruled} remains open.
If $i_1(q)=3$, then $q$ is either a Pfister neighbor or a quasi-Pfister
neighbor by Theorem \ref{neighbor}. If $q$ is quasilinear, then it
is ruled by Theorem \ref{ruled}, and so we can assume that
$q$ is a Pfister neighbor. Then $q$ is special by Witt cancellation
\cite[Theorem 8.3]{EKM} (which shows
that all codimension-1 subforms
of a given Pfister form are similar).
Therefore $q$ is ruled if $i_1(q)=3$. We can expect that there
are no 7-dimensional anisotropic forms with $i_1(q)=2$, as
in characteristic not 2,
but that seems to be unknown. 

More generally, Karpenko
found all possible values of the first Witt index for anisotropic
forms of a given dimension in characteristic not 2 \cite{KarpenkoWitt},
\cite[Theorem 75.13]{EKM}. We can hope
that exactly the same results hold for arbitrary anisotropic
quadratic forms in characteristic 2.

% Omit these bibliography lines if there's no bibliography.

\small \sc DPMMS, Wilberforce Road,
Cambridge CB3 0WB, England

b.totaro@dpmms.cam.ac.uk

\end{document}